\documentclass[11pt,twoside]{amsart}

\usepackage{enumerate,graphics,xypic,floatflt}

\usepackage[hyperindex=true,plainpages=false,colorlinks=false,pdfpagelabels]{hyperref}
\usepackage{chrpaul}

\theoremstyle{plain}
\newtheorem{The}{Theorem}
\newtheorem*{The*}{Theorem}
\newtheorem{Pro}[The]{Proposition}
\newtheorem{Lem}[The]{Lemma}
\newtheorem{Cor}[The]{Corollary}

\theoremstyle{definition}
\newtheorem*{Def}{Definition}

\theoremstyle{remark} 
\newtheorem{Rem}[The]{Remark}
\newtheorem{Exa}[The]{Example}
\newtheorem{Remark}[The]{Remark}

\numberwithin{equation}{section}

\DeclareMathOperator{\Will}{\mathcal{W}}
\DeclareMathOperator{\Area}{\mathcal{A}}
\DeclareMathOperator{\Vol}{\mathcal{V}}
\DeclareMathOperator{\Func}{\mathcal{F}}

\newcommand{\tvector}[1]{{\left(\begin{smallmatrix}#1\end{smallmatrix}\right)}}

\begin{document}

\title{Constrained Willmore Surfaces}

\author{Christoph Bohle}
\author{G.~Paul Peters}
\author{Ulrich Pinkall}

\address{Institut f\"ur Mathematik\\ 
Technische Universit{\"a}t Berlin\\
Stra{\ss}e des 17.\ Juni 136\\
10623 Berlin\\
Germany}

\email{bohle@math.tu-berlin.de\\
  peters@math.tu-berlin.de\\ \newline pinkall@math.tu-berlin.de} 




\thanks{All three authors partially supported by DFG SPP 1154.  The original
  publication is available at
  \href{http://dx.doi.org/10.1007/s00526-007-0142-5}{www.springerlink.com} in
  Calc.\ Var.\ Partial Differential Equations \textbf{32} (2008), 263-277.}

\begin{abstract} 
  Constrained Willmore surfaces are conformal immersions of Riemann surfaces
  that are critical points of the Willmore energy $\Will=\int H^2$ under
  compactly supported infinitesimal conformal variations.  Examples include
  all constant mean curvature surfaces in space forms.  In this paper we
  investigate more generally the critical points of arbitrary geometric
  functionals on the space of immersions under the constraint that the
  admissible variations infinitesimally preserve the conformal structure.
  Besides constrained Willmore surfaces we discuss in some detail examples of
  constrained minimal and volume critical surfaces, the critical points of the
  area and enclosed volume functional under the conformal constraint.
\end{abstract}

\maketitle

\section{Introduction}

In 1882 Felix Klein posed the question whether every abstract Riemann surface
$M$ can be realized in Euclidean 3--space.  A positive answer to this question
was given by A.~Garsia \cite{Ga61} and R.~R{\"u}edy \cite{Ru71}, who
proved that every Riemann surface can be conformally embedded into $\R^3$.
This immediately raises the question for the optimal geometric realizations
in $\R^3$ of a given abstract Riemann surface, or, slightly more general, for
the optimal (especially beautiful, symmetric...) conformal immersions of a
Riemann surface into space.

The same question may be raised about optimal
realizations of topological types of surfaces, without the conformal
constraint. In this context, there is a general agreement that an
excellent way to make the notion of ``optimal'' precise is to look for
surfaces that are critical points of the Willmore functional
$\Will(f)=\int H^2 d\sigma$.

This leads to the notion of constrained Willmore surfaces: an
immersion $f\colon M \to \R^3$ of a Riemann surface $M$ is a
\emph{constrained Willmore surface} if it is critical for $\Will$
under compactly supported infinitesimal conformal variations of $f$.
This generalizes the notion of \emph{Willmore surfaces}, the critical
points of $\Will$ under all compactly supported variations.
Constrained Willmore surfaces are a M{\"o}bius invariant class of
surfaces with strong links to the theory of integrable systems.

The concept of constrained Willmore surfaces already appears in
\cite{PS1,W,Ri97,BPP02,S02}.  However, the Euler--Lagrange equation
characterizing constrained Willmore surfaces is nowhere actually
derived from the variational problem, except for the simplest cases of
compact surfaces of genus $g\leq 1$ in~\cite{BPP02}. Examples of
constrained Willmore surfaces are constant mean curvature surfaces
in 3--dimensional space forms, see \cite{PS1} for the Euclidean case
and \cite{Ri97,BPP02} or Section~\ref{sec:constr-willm-surf} below for
general space forms.

In Section~\ref{sec:constrained_variational_problems}, we derive the
Euler--Lagrange equation for conformal immersions of compact Riemann
surfaces that are critical points of geometric functionals under
infinitesimal conformal variations.  In this Euler--Lagrange equation,
a holomorphic quadratic differential appears as a Lagrange multiplier.
In the case of immersions of non--compact surfaces, the same
Euler--Lagrange equation is still a sufficient (but no longer
necessary) condition for being critical under the conformal
constraint.

In Section~\ref{sec:constr-minim-volume}, we discuss examples of
\emph{constrained minimal} and \emph{constrained volume--critical}
surfaces, i.e., immersions which are critical for area or enclosed
volume under the conformal constraint.  This yields valuable insights
into variational problems of the type studied in
Section~\ref{sec:constrained_variational_problems}.  In particular,
Example~\ref{exa:constrained_no_q} shows that the case of open
surfaces is more delicate than the compact case: we construct a
constrained minimal immersion that does not satisfy the Euler-Lagrange
equation derived in
Section~\ref{sec:constrained_variational_problems}.

In Section~\ref{sec:constr-willm-surf}, we specialize the previous
results to constrained Willmore surfaces and give a short proof of the
fact that constant mean curvature surfaces are constrained Willmore.
Moreover, we describe a phenomenon first observed in~\cite{BoPe2}:
there are (globally smooth) immersions of compact surfaces that are
constrained Willmore after removing a finite number of points, but not
constrained Willmore as compact surfaces.  This phenomenon seems to be
typical for variational problems with conformal constraint:
Example~\ref{exa:discs_of_revolution} in
Section~\ref{sec:constr-minim-volume} shows that the same can be
observed in the case of constrained minimal surfaces.

In the appendix we show that for non strongly isothermic immersions of
compact surfaces all infinitesimal conformal variations may be
realized as derivatives of genuine conformal variations, where
\emph{strongly isothermic} means that the immersion admits a
non--trivial, holomorphic quadratic differential with the property
that the null directions of its real part are principal curvature
directions of the immersion.  As a consequence, a non strongly
isothermic immersion of a compact surface is a critical point of a
geometric functional under infinitesimal conformal variations if and
only if it is critical under genuine conformal variations.

\newpage

\section{Critical Points of Variational Problems\\
  with Conformal Constraint}\label{sec:constrained_variational_problems}

In this section we derive the Euler--Lagrange equation characterizing
immersions of a compact Riemann surface that are critical points of a
geometric functional under the constraint that the conformal structure
is infinitesimally preserved.  In the case of non--compact surfaces,
the equation thus obtained is still a sufficient condition for
critical points, but (as shown in
Section~\ref{sec:constr-minim-volume}) it is not necessary.

Let $M$ be a Riemann surface (without boundary throughout) and denote by
$\Func$ a functional (for non--compact surfaces it is sufficient that the
functional $\Func$ is well defined on compact subsets of $M$) on the space of
immersions of $M$ into some 3--dimensional Riemannian manifold $(\bar M, \bar
g)$. We suppose that $\Func$ is invariant under diffeomorphisms of $M$, i.e.,
that it depends only on the geometry of the immersion and not on the
parametrization.  Moreover, we suppose that $\Func$ is differentiable in the
following sense: for every immersion $ f\colon M\to\bar M$ there exists a
2--form $\grad(\Func) \in \Omega^2(M)$ such that the derivative of $\Func$ in
direction of a compactly supported infinitesimal variation
\begin{equation}
  \label{eq:infinitesimal_def}
  \dot f = u \xi + df(X)\in \Gamma_0(\perp_f\!M \oplus T_f M)
\end{equation}
with $u\in C^\infty_0(M)$ a compactly supported function, $X\in \Gamma_0(TM)$
a compactly supported vector field, and $\xi$ the positive unit normal, is
given by
\[ 
\dot {\!\Func} = \int_M \grad(\Func)u.
\] 
Note that $\,\dot{\! \Func}$ depends only on the normal variation,
because $\Func$ is invariant under diffeomorphisms.

In order to derive the Euler--Lagrange equation for variational problems with
conformal constraint, we need to investigate the effect of
compactly supported infinitesimal variations on the conformal structure. 
For this purpose we use the fact that on oriented surfaces a conformal
structure is the same as a complex structure and can be represented by an
endomorphism field $J\in \Gamma(\End(TM))$ with $J^2=-\Id$.

\begin{Lem}\label{lem:change_of_J} 
  The change of the conformal structure induced by a compactly supported
  infinitesimal normal variation $\dot f = u \xi\in \Gamma_0(\perp_f\! M)$ is
  \begin{equation}
    \label{eq:change_conf}
    \dot J = 2 u \mathring A J,
  \end{equation}
  where $\mathring A$ denotes the trace free part of the Weingarten operator
  $A$.  For a compactly supported infinitesimal tangential variation $\dot f =
  df(X)\in \Gamma_0(T_f\! M)$ the conformal structure changes according to
  \begin{equation}
    \label{change_conf_tang}
    \dot J = \mathcal{L}_X J,
  \end{equation}
  where $\mathcal{L}$ denotes the Lie derivative.
\end{Lem}
\begin{proof}
  Equation \eqref{change_conf_tang} is the definition of the Lie derivative.
  A normal variation $\dot f = u \xi$ changes the induced metric $g$ according
  to
  \[ \dot g = -2 u \,g(A \underline{\; }, \underline{\; }). \] For arbitrary
  vector fields $Y\in \Gamma(TM)$ we have $g(JY,JY)=g(Y,Y)$ and $g(J Y,Y)=0$.
  Infinitesimally, this becomes
  \begin{gather*}
    \dot g(JY,JY) + 2g(\dot JY,JY) = \dot g (Y,Y)\\
    \dot g (JY,Y) + g(\dot JY,Y) = 0.
  \end{gather*}
  Using the equation for $\dot g$ and $\mathring A=\frac12(A+JAJ)$, this
  implies \eqref{eq:change_conf}.
\end{proof}

A compactly supported infinitesimal variation $\dot f = u \xi + df(X)$ is
called \emph{conformal} if the infinitesimal change of the conformal structure
induced by $\dot f$ vanishes. By Lemma~\ref{lem:change_of_J} this is
equivalent to
\begin{equation}
  \label{eq:conformality}
   2 u \mathring A J +\mathcal{L}_X J = 0.
\end{equation}
Clearly, compactly supported infinitesimal conformal variations can be
obtained by taking the derivative $\dot f = {\frac d{dt}}_{|t=0} f_t$ of a
genuine compactly supported conformal variation $f_t$, i.e., a family
$f_t\colon M\to \bar M$ of immersions depending smoothly on $t\in I\subset \R$
such that on the complement of a compact subset of $M$ all $f_t$ coincide with
$f$ and $f=f_0$.  It should be noted that, in contrast to the case of
arbitrary compactly supported infinitesimal variations, it is not true that
every compactly supported infinitesimal conformal variation is the derivative
of a genuine conformal variation: an example of an infinitesimal conformal
variation that does not admit an extension to a genuine conformal variation is
the infinitesimal variation of a planar open annulus $f$ given by a rotational
symmetric function $u\geq 0$ with compact support.  However, in the appendix
we prove that for immersions of compact surfaces that are not strongly
isothermic, every infinitesimal conformal variation is obtained from a genuine
conformal variation.

\begin{Def}
  A conformal immersion $f\colon M \to \bar M$ of a Riemann surface into a
  3--dimensional Riemannian manifold $(\bar M,\bar g)$ is called
  \emph{constrained $\Func$--critical} if
  \[ \dot {\! \Func} = \int_M \grad(\Func)u=0\] for all infinitesimal
  conformal variations $\dot f = u\xi + df(X)$ of $f$ with compact support.
\end{Def}

Accordingly, we call an immersion  a \emph{constrained Willmore surface},
\emph{constrained minimal surface} or  \emph{constrained volume--critical
  surface} if it is constrained $\Func$--critical for $\Func$ the Willmore
functional $\Will$, the area functional $\Area$, or the enclosed volume
functional $\Vol$.

\newpage
For our derivation of the Euler--Lagrange equation we need the adjoint
operators to the operators assigning to a normal or tangential variation
the change of the conformal structure.
We first define the adjoint of the operator $\delta\colon C^\infty_0(M) \to
\Gamma_0(\End_-(TM))$, where $\End_-(TM)$ denotes the bundle of $J$
anti--commuting endomorphisms and where
\begin{equation}\label{eq:delta}
  \delta(u)=2u\mathring A J
\end{equation} 
is the change \eqref{eq:change_conf} of the complex structure induced by a
normal variation $u\in C^\infty_0(M)\cong \Gamma_0(\perp_f\!M)$.  For this we
need the non--degenerate pairing 
\begin{equation}\label{eq:pairing1}
  \<\omega,u\> = \int_M \omega u 
\end{equation}
between 2--forms $\Omega^2(M)$ and normal variations $C^\infty_0(M)\cong
\Gamma_0(\perp_f\!M)$ and the non--degenerate pairing
\begin{equation}
  \label{eq:pairing2}
  \<q,R\> =\int_M 2\Re(q)(R\underline{\;}\land \underline{\; })
\end{equation}
between the space of quadratic differentials $\Gamma(K^2)$ and
$\Gamma_0(\End_-(M))$, where we use the convention that, for $b$ a bilinear
form, the expression $b(\underline{\;}\land \underline{\; })$ denotes the
2--form $\omega$ defined by $\omega(X,Y) = b(X,Y)-b(Y,X)$. With respect to
these pairings the adjoint operator $\delta^*\colon \Gamma(K^2) \to
\Omega^2(M)$ of $\delta$ is given by
\begin{equation}\label{eq:delta_star}
  \delta^*(q) = 4 \Re(q)(\mathring A J\underline{\;}\land \underline{\; }). 
\end{equation}
The spaces involved in the definition of $\delta^*$ are visualized by the
diagram
\begin{equation}
  \label{eq:diag1} 
  \xymatrix{{\qquad C^\infty_0(M)\qquad} \ar[r]^\delta \ar@{<.>}[d] &
    \Gamma_0(\End_-(TM))\ar@{<.>}[d] \\ 
    {\qquad \Omega^2(M) \qquad}  & {\qquad \Gamma(K^2),\qquad} \ar[l]^{\delta^*} }
\end{equation}
where the vertical arrows indicate the pairings defined in \eqref{eq:pairing1}
and \eqref{eq:pairing2}.

The operator $X\mapsto \mathcal{L}_X J$ assigning to a tangential variation
the induced change of the complex structure can be interpreted as follows:
under the canonical isomorphisms $TM\cong K^{-1}$ and $\End_-(TM) \cong \bar K
K^{-1}$, with $K$ the canonical bundle of the Riemann surface $M$, the
operator $X\mapsto \mathcal{L}_X J$ is the usual $\dbar$--operator on vector
fields. We therefore denote it by
\begin{equation}
  \label{eq:dbar1}
  \dbar \colon \Gamma_0(TM)\cong\Gamma_0(K^{-1})\to
  \Gamma_0(\End_-(TM))\cong \Gamma_0(\bar K K^{-1}).
\end{equation}
Its dual operator is as well a $\dbar$--operator which we denote by 
\begin{equation}
  \label{eq:dbar2}
  \dbar^* \colon \Gamma(K^2) \to \Gamma(\bar K K^2) \cong
  \Omega^2(T^*M).  
\end{equation}
It coincides with the usual $\dbar$--operator on quadratic differentials.  As
in the case of normal variations, we visualize this by the diagram
\begin{equation}
  \label{eq:diag2} 
  \xymatrix{ {\qquad \Gamma_0(TM) \qquad} \ar[r]^\dbar \ar@{<.>}[d] &
    \Gamma_0(\End_-(TM))\ar@{<.>}[d] \\ 
    {\qquad \Omega^2(T^*M)\qquad}  & {\qquad
      \Gamma(K^2). \qquad}\ar[l]^{\dbar^*} }  
\end{equation}

As an application of this language we obtain the following sufficient
condition for constrained $\Func$--critical immersions.
 
\begin{Pro}\label{pro:sufficient}
  Let $f\colon M \to \bar M$ be a conformal immersion of a Riemann surface
  into a 3--dimensional Riemannian manifold $(\bar M,\bar g)$.  If there is a
  holomorphic quadratic differential $q\in H^0(K^2)$ such that $\grad(\Func) =
  \delta^*(q)$, then $f$ is a constrained $\Func$--critical immersion.
\end{Pro}
\begin{proof}
  Let $\dot f= u \xi+df(X) \in \Gamma_0(\perp_f\! M\oplus T_fM)$ be a
  compactly supported infinitesimal conformal variation. By
  \eqref{eq:conformality}, $\delta(u) + \dbar(X) =0$ and therefore
  \[ \dot {\!\Func} = \< \grad(\Func), u \> = \< \delta^*(q), u \> = \< q,
  \delta(u)\> =0,
  \] 
  where the last equality follows from $\delta(u)\in \im(\dbar)\subseteq
  \ker(\dbar^*)^\perp$.
\end{proof}

For the case of compact surfaces the following theorem shows that this
sufficient condition is also necessary.  The theorem follows essentially from
the fact that, by elliptic theory of the $\dbar$--operator (or the so called
Weyl--Lemma), for compact surfaces the space $H^0(K^2)=\ker(\dbar^*)$ of
holomorphic quadratic differentials is finite dimensional and $\im(\dbar) =
\ker(\dbar^*)^\perp$.

\begin{The}
  Let $f\colon M \to \bar M$ be a conformal immersion of a compact
  Riemann surface.  Then $f$ is constrained $\Func$--critical if and only if
  there is a holomorphic quadratic differential $q\in H^0(K^2)$ such that
  \[ \grad(\Func) = \delta^*(q).\]
\label{Th:Euler_Lagrange}\end{The}
\begin{proof}
  By definition, the immersion $f$ is constrained $\Func$--critical if and
  only if $\<\grad(\Func),u\>=0$ for all normal deformations $u\xi \in
  \Gamma(\perp_f\! M)$ such that $\delta(u)\in \im(\dbar)$. From $\im(\dbar) =
  \ker(\dbar^*)^\perp$ we obtain that $\delta(u)\in \im(\dbar)$ is equivalent
  to $u \in (\delta^*(\ker(\dbar^*)))^\perp$.  

  Hence, $f$ is constrained $\Func$--critical if and only if $\grad(\Func) \in
  (\delta^*(\ker(\dbar^*)))^{\perp\perp}$. This proves the theorem, because
  $\ker(\dbar^*)$ is finite dimensional and therefore
  $(\delta^*(\ker(\dbar^*)))^{\perp\perp} = \delta^*(\ker(\dbar^*))$.
\end{proof}

\begin{Remark}
  We refer to $q$ as the \emph{Lagrange multiplier} of the constrained
  $\Func$--critical immersion, because $\grad(\Func) = \delta^*(q)$ can be
  interpreted as the condition that $\grad(\Func)$ is orthogonal to the
  space of conformal immersions $f\colon M\to\bar M$, see the appendix for
  details.
\end{Remark}

We conclude the section by briefly comparing three different notions
of critical points of a geometric functional $\mathcal{F}$ on the
space of immersions under the constraint that the conformal structure
is fixed:
\begin{enumerate}
\item[i)] The immersions satisfies the Euler--Lagrange equation 
  \[\grad(\Func) =  \delta^*(q)\]
  for some Lagrange multiplier $q\in H^0(K^2)$.
\item[ii)] The immersion is constrained $\Func$--critical, i.e., it is critical
  under all infinitesimal conformal variations with compact support.
\item[iii)] The immersion is  $\Func$--critical under all genuine conformal variations with
  compact support.
\end{enumerate}
Proposition~\ref{pro:sufficient} shows that i) always implies ii) and,
by Theorem~\ref{Th:Euler_Lagrange}, if the underlying surface is compact
i) and ii) are equivalent (we will see in
Example~\ref{exa:constrained_no_q} below that this is
not the case for non--compact surfaces).  By definition, ii) always
implies iii).  Corollary~\ref{cor:non_isothermic} in the appendix
shows that ii) and iii) are equivalent for immersions of compact
surfaces that are not strongly isothermic.

\section{Constrained Minimal and Volume--Critical Surfaces}\label{sec:constr-minim-volume}

In this section we discuss examples of \emph{constrained minimal surfaces} and
\emph{constrained volume--critical surfaces}, the critical points of
area~$\Area$ and enclosed volume~$\Vol$ under compactly supported
infinitesimal conformal variations.  These examples reveal several remarkable
properties of solutions to variational problems under the constraint that the
admissible variations preserve the conformal structure:
\begin{itemize}
\item Example~\ref{exa:spheres_of_revolution} shows that there are
  immersions of compact surfaces  that are not constrained
  $\mathcal{F}$--critical, but become constrained
  $\mathcal{F}$--critical after removing a finite number of points.
\item Example~\ref{exa:cylinder} and \ref{exa:hopf_cylinders} show
  that constrained minimal immersions are not necessarily analytic.
  This reflects the fact that the Lagrange multiplier appearing
  in the Euler--Lagrange equation for constrained minimal surfaces
  destroys its ellipticity.
\item Example~\ref{exa:constrained_no_q} shows that, if the underlying
  Riemann surface is non--compact, constrained $\mathcal{F}$--critical
  immersions do not necessarily admit a Lagrange multiplier $q\in
  H^0(K^2)$ such that $\grad(\mathcal{F}) = \delta^*(q)$.
\end{itemize}

\begin{Lem}\label{L:area_volume_gradient}
  Let $f\colon M\to\bar M$ be an immersion of a Riemann surface $M$ into a
  3--dimensional Riemannian manifold $(\bar M,\bar g)$. Under the
  compactly supported infinitesimal normal variation $\dot f = u \xi$,
  the area and volume functional change according to
  \begin{align*}
    \dot{\!\! \Area}  = - 2\int_M u H \,d\sigma \qquad\quad \textrm{and} 
    \qquad\quad   \dot \Vol  =  \int_M u \,d\sigma,
  \end{align*}
  where $H$ denotes the mean curvature and $d\sigma$ the area element of the
  immersion.  Equivalently, the gradients of the functionals are the 2--forms
  \begin{align*}
    \grad(\Area) = -2 H\, d\sigma  \qquad\quad \textrm{and} 
    \qquad\quad \grad(\Vol) = d\sigma.
  \end{align*}
\end{Lem}

A proof of this lemma is left to the reader. When considering the
functionals $\Area$ and $\Vol$ without the conformal constraint they
behave rather differently: while $\Vol$ has no critical points at all,
the critical points of the functional $\Area$ are the minimal
immersions into $(\bar M,\bar g)$.

\begin{Remark}
  In the context of the isoperimetric problem one is interested in the
  critical points of $\Area$ under the constraint that $\Vol$ is
  fixed
  (or vice versa).  Solutions to this constrained variational problem are the
  surfaces of constant mean curvature, since $f$ is critical under the
  constraint if and only if $\grad(\Area) =-2 H \, d\sigma \in \Span\{ d\sigma
  \}^{\perp\perp}$.
\end{Remark}

In Euclidean space there exist no immersions of compact Riemann surfaces that
are constrained minimal or constrained volume--critical, because homotheties
do not change the conformal structure but scale $\Area$ and $\Vol$. In
contrast to this, there are compact constrained minimal and constrained
volume--critical surfaces in the 3--sphere.
\begin{Exa}[Homogeneous Surfaces in Space Forms]
  \emph{Every non totally umbilic homogeneous surface in a space form of
    dimension $3$ is constrained minimal and constrained volume--critical.}
  This can be seen as follows.  Because homogeneous surfaces have constant
  mean curvature, the Codazzi equation implies that the Hopf differential,
  i.e., the unique quadratic differential $Q\in \Gamma(K^2)$ that satisfies
  $\mathring \sff=(Q+\bar Q)\xi$, where $\mathring \sff$ denotes the trace free
  second fundamental form, is holomorphic. From $\Re(Q)=\tfrac12 g(\mathring A
  \underline{\;}, \underline{\;})$, $\mathring A^2 = (H^2 -G)\Id$ and
  $d\sigma=\tfrac12 g(J\underline{\;}\land \underline{\; })$ follows that the
  Hopf differential satisfies
  \begin{equation}
    \label{eq:hopf_in_delta}
    \delta^*(Q)=4(H^2-G) d\sigma,
  \end{equation}
  where $G$ denotes the Gaussian curvature, i.e., $G=\det(A)$.  Since $H^2-G$
  is a nonzero constant, there are $\lambda_1$ and $\lambda_2\in \R$ such
  that $\grad(\Area)=\lambda_1\delta^*(Q)$ and
  $\grad(\Vol)=\lambda_2\delta^*(Q)$, which, by
  Proposition~\ref{pro:sufficient}, proves the statement.
  
  The non totally umbilic homogeneous surfaces in space forms are isometric
  to the following examples: the tori $(u,v)\mapsto (r_1 e^{{\i} u}, r_2
  e^{{\i} v})$ with $r_1^2+r_2^2 = 1$ in $S^3\subset \C^2$, the cylinders
  $(u,v)\mapsto (r\cos(u),r\sin(u), v)$ in $\R^3$ and the cones in the half
  space model of hyperbolic 3--space whose cusp lies on the infinity boundary
  and whose axis of rotation is perpendicular to the infinity boundary.
  
  The exclusion of totally umbilic surfaces in the above statement is
  necessary: every piece of the plane is minimal but not constrained
  volume--critical, because every infinitesimal deformation with compact
  support is conformal. For the same reason, every piece of the round sphere
  is neither constrained minimal nor constrained volume--critical.
\end{Exa}

\begin{Exa}[Discs of revolution]\label{exa:discs_of_revolution}
  Let $f\colon \Delta\to\R^3$ be a conformal immersion of the open unit disc
  into $\R^3$ that, on the punctured disc $\Delta\setminus\{0\}$, parametrizes
  a surface of re\-vo\-lu\-tion. Then every compactly supported variation of
  $f$ that preserves the rotational symmetry is conformal.

  \begin{floatingfigure}{4cm}\centering\ \\[2ex] 
    \scalebox{.6}{\includegraphics{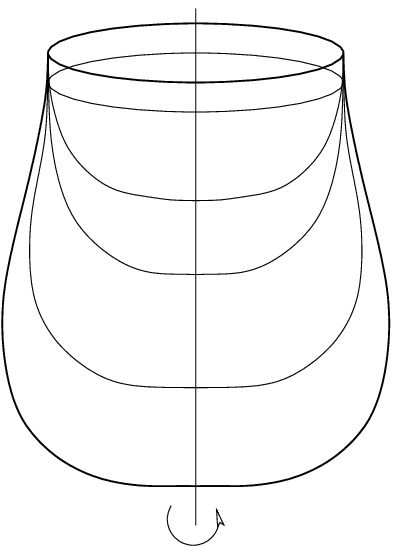}}\\
    {\small \textsc{Figure 1.} Deformation of Surfaces of
      Revolution.} \\[1.5ex] 
  \end{floatingfigure}
  The reason is that for immersions of the disc with rotational
  symmetry, the meridian curve has infinite length in the hyperbolic
  half plane whose boundary is the axis of rotation, and, moreover,
  such immersions are conformal if and only if, with respect to polar
  coordinates $z=e^{x+{\i} y}$, the meridian curve is pa\-ra\-me\-trized
  by hyperbolic arc length.  Thus, deforming the profile curve one can
  decrease $\Area$ and $\Vol$ and therefore: \emph{An immersion of a
    disc with rotational symmetry is never constrained minimal nor
    constrained volume--critical, unless it is a
    planar disc and therefore minimal.}
\end{Exa}

\begin{Exa}[Spheres of revolution]\label{exa:spheres_of_revolution}
Figure~2 shows four constrained minimal surfaces with rotational symmetry
\cite{Sz04}. These surfaces can be smoothly extended to compact embedded
spheres which are not constrained minimal. The preceding argument shows that
only after removing both points on the axis of revolution one obtains
  constrained minimal surfaces.

\begin{center}
\scalebox{.4}{\includegraphics{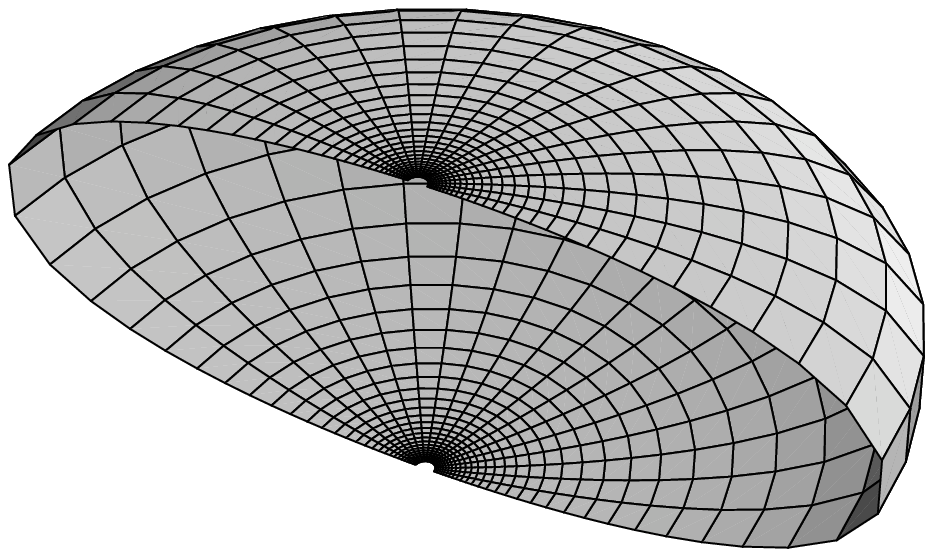}}\hfill 
\scalebox{.35}{\includegraphics{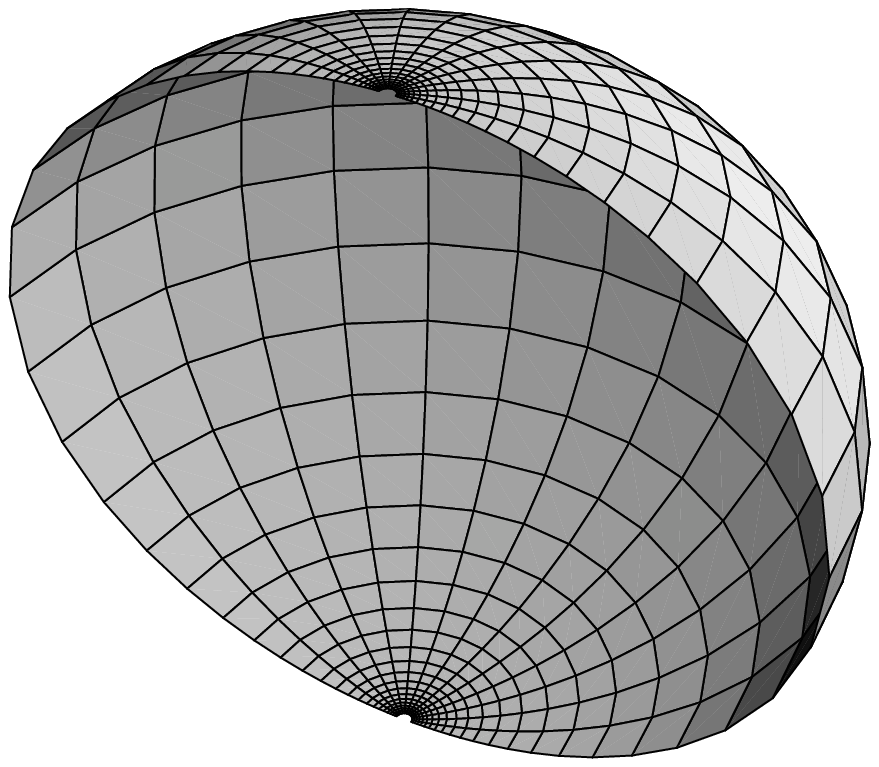}}\hfill 
\scalebox{.35}{\includegraphics{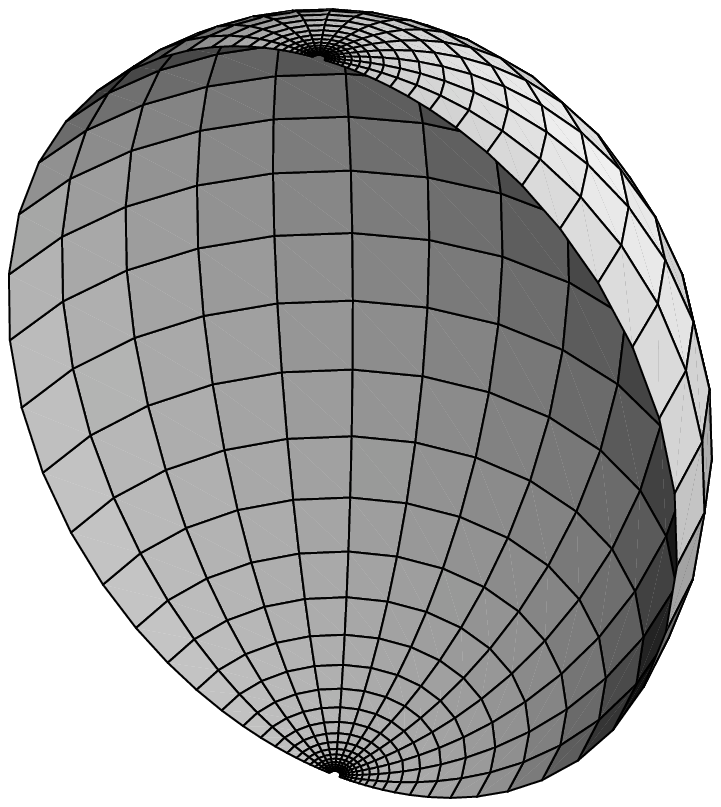}}\hfill 
\scalebox{.4}{\includegraphics{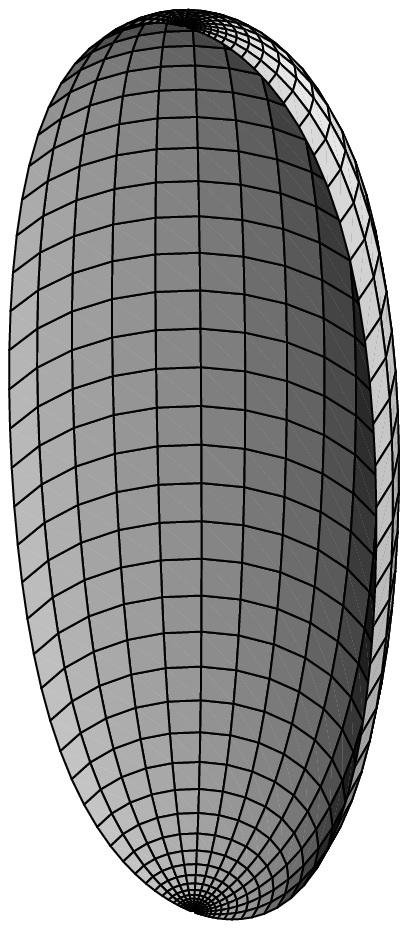}}

{\small \textsc{Figure 2.} Constrained minimal surfaces.}
\end{center}

Interestingly, the round sphere with poles removed can be approximated
by such punctured constrained minimal spheres with rotational symmetry
although, by Example~\ref{exa:discs_of_revolution}, a piece of the
sphere itself is never constrained minimal.
\end{Exa}

\begin{Exa}[Cylinders over plane curves]\label{exa:cylinder}
  This example shows that, in contrast to minimal surfaces, constrained
  minimal surfaces are not necessarily analytic. \emph{Every cylinder over a
    plane curve is constrained minimal.\/} To see this, assume that the plane
  curve $\gamma(t)=(x(t),y(t))$ is parametrized with respect to arc length.
  Then the cylinder $f(u,v) = (x(u),y(u),v)$ is an isometric immersions of the
  plane with Weingarten operator $A= \tvector{ -\kappa & 0 \\ 0 & 0}$ and mean
  curvature $H=-\tfrac12 \kappa$ where $\kappa$ denotes the curvature of
  $\gamma$. For the holomorphic quadratic differential $q=dz^2$ with $z=u+{\i}
  v$ we have $\delta^*(q)= -4\kappa \,d\sigma$ and therefore
  \begin{equation*}
      \grad(\Area) = -2 H\, d\sigma = \kappa \,d\sigma = -\tfrac14 \delta^*(q)
  \end{equation*}
  which, by Proposition~\ref{pro:sufficient}, proves the statement.
\end{Exa}

\begin{Exa}[Constrained minimal surface without Lagrange multiplier
  $q$]\label{exa:constrained_no_q} 
  We give now an example of a constrained minimal immersion that does
  not admit a holomorphic quadratic differential $q$ with
  $\grad(\Area)=\delta^*(q)$.  This shows that
  Theorem~\ref{Th:Euler_Lagrange} does not hold in the non--compact
  case. The example is constructed as follows: take a planar domain as
  in Figure 3 and bend the hatched areas upwards such that all 3
  fingers of the surface are cylinders over plane curves lying in
  three different planes perpendicular to the article.

  As shown in Example~\ref{exa:cylinder} the surface thus obtained is locally
  constrained minimal and locally admits a holomorphic quadratic differential
  $q$ such that the Euler--Lagrange equation $\grad(\Area)=\delta^*(q)$ is
  satisfied.  It remains to be proven that the surface is globally constrained
  minimal but does not admit a global holomorphic quadratic differential $q$
  such that $\grad(\Area)=\delta^*(q)$.

   \begin{floatingfigure}{4.4cm}\centering
    \scalebox{.4}{\includegraphics{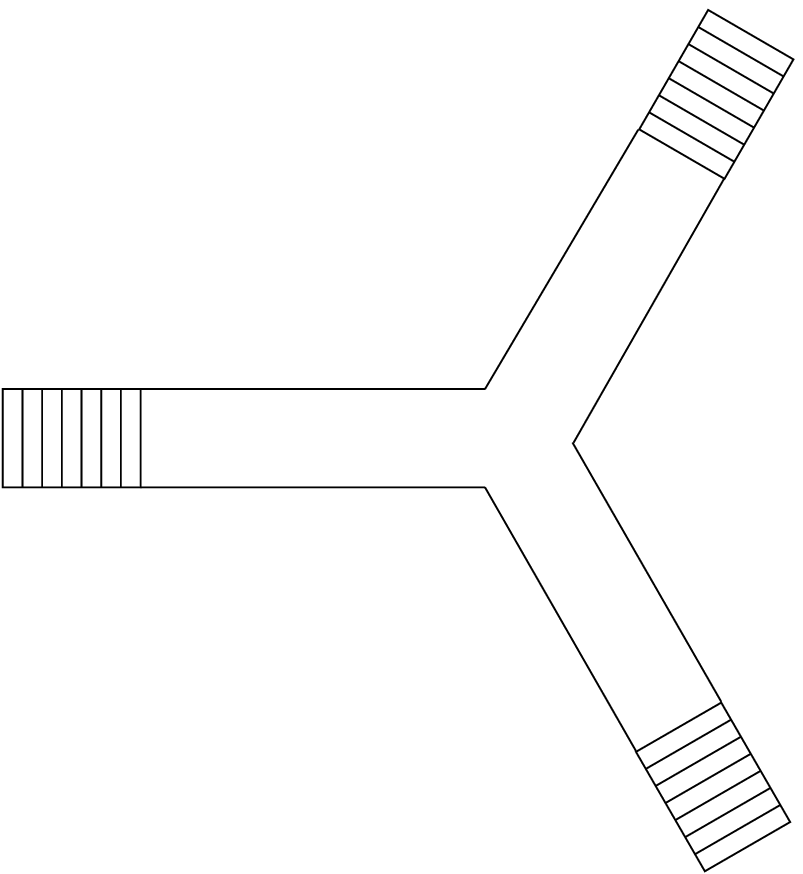}}
    {\small \textsc{Figure 3.} Constrained minimal surface without $q$.}
  \end{floatingfigure}
 To prove the first claim, we have to check that $\int_M \grad(\Area) u=0$
  for every compactly supported infinitesimal conformal variation $\dot f =
  u\xi + df(X)$.  On the white planar area $\grad(\Area)$ vanishes and
  $\delta(\tilde u)$ vanishes for all normal variations $\tilde u$.  Hence,
  without loss of generality we may assume that $u$ has support in a small
  neighborhood of the hatched areas: changing $u$ on the white domain to
  achieve this would neither change the integral $\int_M \grad(\Area) u$ nor
  would it destroy conformality of the variation.  But this proves $\int_M
  \grad(\Area) u=0$, because the 3 fingers of the surface are constrained
  minimal by
  Example~\ref{exa:cylinder}. \\
  \phantom{\quad} The fact that there is no global holomorphic quadratic
  differential $q$ such that $\grad(\Area)=\delta^*(q)$ can be seen as
  follows: as in Example~\ref{exa:cylinder}, denote by $q_j$, $j=1,...,3$ the
  squares of the differentials of the cylinder coordinates on the three
  fingers. On each finger we then have $\grad(\Area)=-\frac14\delta^*(q_j)$
  for the respective $j$. Moreover, if there was a global holomorphic
  quadratic differential $q$ with $\grad(\Area)=\delta^*(q)$, then on each
  finger there had to be $t_j\in\R$ such that $-4q=q_j+t_jiq_j$, which is
  impossible.
\end{Exa}

\begin{Exa}[Hopf cylinders in $S^3$]\label{exa:hopf_cylinders}
  \emph{All Hopf cylinders in $S^3$, in particular all Hopf tori, are
    constrained minimal and constrained volume--critical}.  A surface in $S^3$
  that is the preimage $\pi^{-1}(\gamma)$ of an immersed curve $\gamma$ in
  $S^2$ under the Hopf fibration $\pi\colon S^3\to S^2$ is called a Hopf
  cylinder or, if $\gamma$ is closed, a Hopf torus.  Taking arc length
  parameters of a horizontal lift of $\gamma$ and of the Hopf fibers yields an
  isometric parametrization of the Hopf cylinder.  With respect to the
  corresponding coordinates $z=x+i y$, the Weingarten operator takes the form
  $A=\tvector{-2\kappa & -1\\-1&0}$, cf.\ \cite{P85}, and the mean curvature
  $H=-\kappa$, where $\kappa$ denotes the curvature of $\gamma$ in $S^2$.  The
  holomorphic quadratic differential $q=dz^2$ satisfies $\delta^*(q)=-8\kappa
  d\sigma$ and $\delta^*(iq)=8d\sigma$. By Lemma~\ref{L:area_volume_gradient}
  and Proposition~\ref{pro:sufficient} this proves the claim.

    Using Theorem~\ref{Th:constrained_Willmore} of the next section one gets
  that \emph{a Hopf torus is constrained Willmore if and only if there are constants
  $a,b\in\R$ such that}
  \begin{equation*}
    2\kappa''+\kappa^3 +a \kappa +b=0.  
  \end{equation*}
  This equation means that the generating closed curve $\gamma$ is a
  critical point of $\int_\gamma\kappa^2 ds$ for variations fixing the length and
  the enclosed area, i.e., \emph{$\gamma$ is a generalized elastic curve},
  cf.~\cite{LS84} .  This reflects the fact, see \cite{P85}, that the
  Hopf torus generated by $\gamma$ has Willmore energy $\Will=\pi \int_\gamma
  (\kappa^2 + 1)ds$ and that fixing the conformal type of a Hopf torus
  means fixing the length $L$ and enclosed area $A$ of $\gamma$, because
  the Hopf torus is isometric to $\R^2/\Gamma$ with $\Gamma$ generated by
  $(2\pi,0)$ and $(A/2,L/2)$.
\end{Exa}

\section{Constrained Willmore Surfaces}
\label{sec:constr-willm-surf}

The Willmore energy or elastic bending energy was suggested as a global
invariant of surfaces by T.~Willmore~\cite{W0} in 1965, but appears already in
earlier work of S.~Germain~\cite{G} and W.~Blaschke~\cite{B}.  Its
applications range from the biophysics of membranes to string theory, where it
was introduced by W.~Helfrich and A.~M.~Polyakov, respectively. The Willmore
energy of an immersion $f\colon M \to \bar M$ into a 3--dimensional Riemannian manifold
$(\bar M,\bar g)$ is
\begin{equation}
  \label{eq:will_general_def}
  \Will(f) = \int_M (H^2+\bar K) d\sigma,
\end{equation}
where $H$ is the mean curvature and $\bar K$ the sectional curvature of $\bar
M$ along $T_fM$. Using the Gauss equation $K= G + \bar K$ which relates the
sectional curvature $K$ of the induced metric on $M$ to the Gaussian curvature
$G=\det(A)$ of the immersion the Willmore functional can be rewritten as
\begin{equation}
  \label{eq:will_moebiusinv}
  \Will(f) = \underbrace{\int_M |\mathring\sff|^2
    d\sigma}_{=: \tilde \Will} + \int_M K d\sigma,
\end{equation}
where $|\mathring\sff|^2= (H^2 -G)$ denotes the squared length of the trace
free second fundamental form.  The functional $\tilde \Will$ is invariant
under conformal transformations, because the integrand $|\mathring \sff|^2
d\sigma$ does not change under conformal scalings of the ambient metric. By
the Gauss--Bonnet theorem, since we are considering only compactly supported
variations, being a critical point of the Willmore functional is a conformally
invariant property, and, for compact surfaces, the Willmore functional $\Will$
itself is conformally invariant.

\begin{Remark} The Willmore functional $\Will$ and the functionals $\int_M
  (\kappa_1^2+\kappa_2^2)d\sigma$ and $\int_M (\kappa_1-\kappa_2)^2d\sigma$,
  with $\kappa_1$, $\kappa_2$ denoting the principal curvatures, differ by
  multiples of $\int_M K d\sigma$ only and therefore have the same critical
  points.  Thus, the Willmore functional can be seen as a measure of the total
  amount of principal curvature and the defect from being totally umbilic.
\end{Remark}

\begin{Def}
  A conformal immersion $f\colon M\to \bar M$ of a Riemann surface $M$ into a
  3--dimensional Riemannian manifold $\bar M$ is \emph{constrained Willmore}
  if it is a critical point of the Willmore functional $\Will$ under compactly
  supported infinitesimal conformal variations.
\end{Def}

For simplicity we restrict now to the case that $(\bar M,\bar g)$ is a
manifold of constant sectional curvature.  Because both the functional
and the constraint are invariant under conformal changes of the metric
$\bar g$, the notion of constrained Willmore surfaces for all three
space forms coincides and could as well be considered from a purely
M{\"o}bius geometric viewpoint, see e.g.~\cite{BPP02}.

A proof of the following theorem can be found in \cite{We78}.
\begin{The}\label{Th:gradient_will}
  A compactly supported infinitesimal variation $\dot f = u \xi + df(X)$ of an
  immersion $f\colon M \to \bar M$ into a 3--dimensional space form changes
  the Willmore functional according to
  \begin{equation}
    \label{eq:variation_willmore}
    \dot \Will = \int_M u(\Delta H + 2H(H^2-G))d\sigma.
  \end{equation}
  In particular, the gradient of the Willmore functional is the 2--form
\begin{equation}
  \label{eq:gradient_willmore}
  \grad(\Will) = (\Delta H + 2H(H^2-G))d\sigma.
\end{equation}
\end{The}

Proposition~\ref{pro:sufficient} and Theorems~\ref{Th:Euler_Lagrange}
and~\ref{Th:gradient_will} imply the following theorem.

\begin{The}\label{Th:constrained_Willmore}
  A conformal immersion of a compact Riemann surface into a 3--dimensional
  space form is constrained Willmore if and only if there exists a holomorphic
  quadratic differential $q\in H^0(K^2)$ such that
\begin{equation}\label{eq:euler-lagrange-cws}
  (\Delta H + 2H(H^2-G))d\sigma = \delta^*(q).
\end{equation}
For non--compact surfaces, the existence of a holomorphic quadratic
differential $q\in H^0(K^2)$ that satisfies \eqref{eq:euler-lagrange-cws}
implies that the surface is constrained Willmore.
\end{The}

It should be noted that Willmore surfaces correspond to the case that
$q=0$ in \eqref{eq:euler-lagrange-cws}. In fact, even for immersions
of non--compact surfaces the Euler--Lagrange equation of Willmore
surfaces is a necessary condition: an immersion of a (possibly open)
surface into a 3--dimensional space form is Willmore if and only if $\Delta
H + 2H(H^2-G)=0$.

We now give a simple proof of the following result, cf.~\cite{PS1,Ri97,BPP02}.

\begin{Cor}\label{Cor:cmc_is_constrained_Willmore}
  Every constant mean curvature surface $f\colon M\to \bar M$ in a
  3--dimensional space form $(\bar M,\bar g)$ is constrained Willmore.
\end{Cor}

\begin{proof}
  For constant mean curvature surfaces, the gradient of $\Will$ is given by
  $\grad(\Will) = 2 H (H^2-G) d\sigma$.  Moreover, by
  \eqref{eq:hopf_in_delta}, the holomorphic Hopf differential $Q$ satisfies
  $\delta^*(Q)=4(H^2-G)d\sigma$. This proves the statement, because the
  holomorphic quadratic differential $q=\tfrac12 H Q$ satisfies \[
  \grad(\Will) = \delta^*(q). \vspace{-0.5cm} \] 
\end{proof}

In addition to being constrained Willmore, constant mean curvature surfaces in
space forms are isothermic, which means that away from umbilics they admit
conformal curvature line coordinates.  Hence, constant mean curvature surfaces
belong to two important surface classes of M{\"o}bius geometry: that of
constrained Willmore surfaces and that of isothermic surfaces.  For a
constrained Willmore immersion of a torus into the conformal 3--sphere,
J.~Richter conversely proved (cf.~\cite{BPP02}) that, if the immersion
is also strongly isothermic as defined in the appendix, it has
constant mean curvature with respect to some space form subgeometry.  This
generalizes a theorem of G.~Thomsen~\cite{B}: a surface is Willmore and
isothermic if and only if it is minimal in some space form.

In contrast to Thomsen's  theorem, the global assumption in Richter's
result is essential: an example, which is due to F. Burstall, of an
isothermic, constrained Willmore surface that does not have constant mean
curvature in some space form is the cylinder over the plane curve\footnote{The
  curvature of this plane curve is the solution to $\kappa''(s) +
  \tfrac12\kappa(s)^3 = (a+b s)\kappa(s)$ with $a=0.2$, $b=0.02$ and
  $\kappa(0)=1$, $\kappa'(0)=0$.} in Figure 4.

\begin{center}
  \scalebox{.6}{\includegraphics{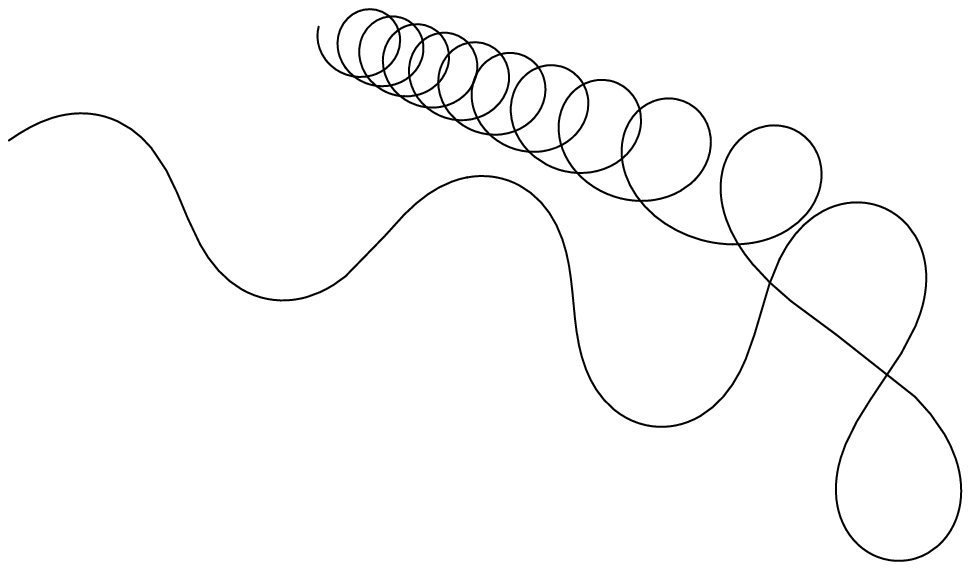}}\\
  {\small \textsc{Figure 4.} Generating curve of a Burstall--cylinder.}
\end{center}

A fundamental property of constrained Willmore surfaces is the
existence of an associated family depending on a spectral parameter,
see e.g.~\cite{BPP02}.  This relation to the theory of integrable
systems has important consequences.  For example, in~\cite{Bo} the
associated family of flat connections is used in order to show that to every
constrained Willmore torus one can assign a compact Riemann surface,
the so called spectral curve (a different approach to proving this is
described in~\cite{S02}).  This fact allows to parameterize constrained
Willmore tori explicitly in terms of holomorphic functions.

The case of constrained Willmore immersions of the sphere is the only case
where a complete classification (in the sense of a reduction to a simpler
algebraic geometric problem) is known: because there is only one conformal
structure on the sphere, all constrained Willmore spheres are Willmore and one
can apply R.~Bryant's result \cite{Br84,Br88} according to which all Willmore
spheres in the conformal 3--sphere are complete Euclidean minimal surfaces of
finite total curvature with planar ends (and therefore algebraic). Bryant
proved that the Willmore energy of Willmore spheres is quantized and that the
possible Willmore energies are $W\in 4\pi (\N^* \backslash\{2,3,5,7\})$. From
the integrable systems point of view it is interesting that Willmore
spheres are examples of so called soliton spheres, see \cite{BoPe1}.

Based on the idea of smooth ends \cite{BoPe2} for constant mean
curvature 1 surfaces in hyperbolic 3--space one can construct a family
of globally smooth conformal immersions of the sphere that, by
Corollary~\ref{Cor:cmc_is_constrained_Willmore}, are constrained
Willmore after removing a finite number of points --- the ``ends'' at which
the immersion touches the ideal boundary of hyperbolic space.  These so
called Bryant spheres with smooth ends are constrained Willmore as
punctured spheres only, but not as compact surfaces, because then they
had to be Willmore which is impossible by Bryant's result
on Willmore spheres.  Nevertheless, they obey the same quantization of the
Willmore energy as Willmore spheres, i.e., $W\in 4\pi (\N^* \backslash\{2,3,5,7\})$,
and they also are soliton spheres.  Figure~5 shows a Bryant sphere
with two smooth ends.

\begin{center}
    \scalebox{.6}{\includegraphics{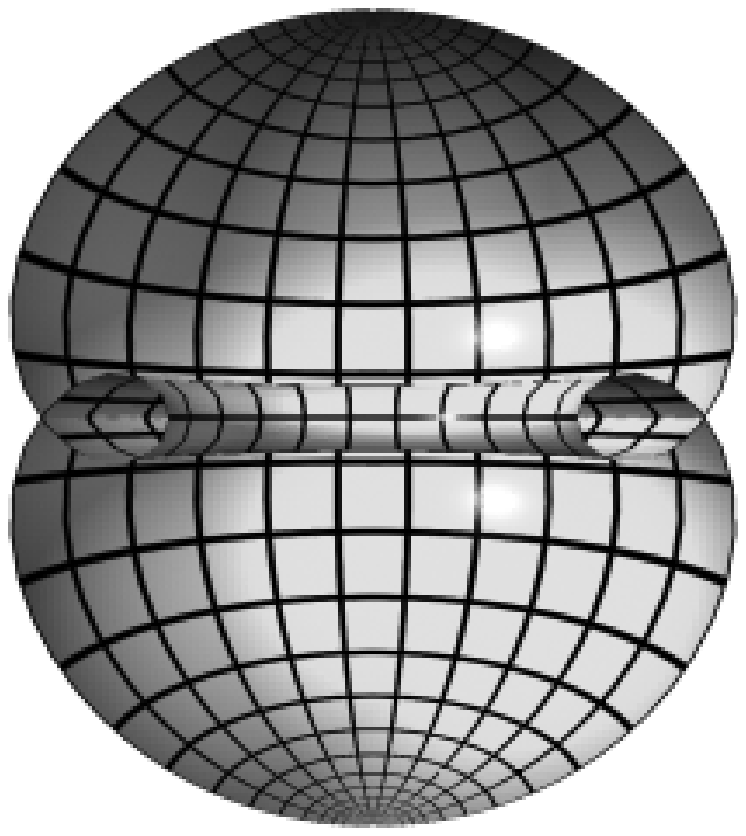}}\\
   {\small \textsc{Figure 5.} Catenoid cousin with smooth ends.}
\end{center}

The phenomenon that an immersion of a surface becomes constrained
Willmore only after removing a finite number of points does not occur
for Willmore surfaces, i.e., for the critical points of the
Willmore energy under all variations. In the example
of Bryant spheres with smooth ends, the holomorphic quadratic
differential in the Euler--Lagrange equation has poles at the ends,
which have to be removed in order to obtain a constrained Willmore
surface.  Because the left hand side in \eqref{eq:euler-lagrange-cws}
is globally smooth, the ends have to be umbilic points.

\appendix
\section*{Appendix. Infinitesimal and genuine conformal transformations}

Our definition of constrained $\Func$--critical immersions in
Section~\ref{sec:constrained_variational_problems} demands that the
immersion is critical for all \emph{infinitesimal} conformal
variations with compact support.  Constrained $\Func$--critical
immersions according to this definition are in particular critical
under all \emph{genuine} conformal variations with compact support.
In this appendix we show that for non strongly isothermic immersions
of compact surfaces it is conversely true that being critical with
respect to genuine variations implies being critical for all infinitesimal
conformal variations.

\begin{Def}
  We call a conformal immersion $f\colon M\to \bar M$ of a Riemann surface
  $M$ into a 3--dimensional Riemannian manifold $\bar M$
  \emph{strongly isothermic} if there exists a non--trivial
  holomorphic quadratic differential $q\in H^0(K^2)$ with the property
  that the null directions of $\Re(q)$ are principal curvature
  directions of $f$.
\end{Def}

An alternative characterization of strongly isothermic surfaces is provided by
the following proposition.

\begin{Pro}\label{pro:isothermic_is_deltastar_of_q_eq_0}
  Let $f$ be a conformal immersion of a Riemann surface and let $q\in
  H^0(K^2)$ be a non--trivial holomorphic quadratic differential.  The
  null directions of $\Re(q)$ are principle curvature directions of $f$
  if and only if \[\delta^*(q)=0\] for $\delta^*$ as in \eqref{eq:delta_star}.
  The zeros of $q$ are then umbilic points of $f$.
\end{Pro}

\begin{proof}
  Away from the isolated zeros of $q$ one can choose local coordinates
  $z=x+iy\colon M\supset U\to\C$ such that $dz^2=iq$ on $U$.  Then $\Re(q)=2dx dy$
  and we have
  \begin{equation*}
    \delta^*(q)= 8 F\, dx\land dy,
  \end{equation*}
  where $F$ denotes the off--diagonal term of the Weingarten operator
  $A$ with respect to the coordinates $(x,y)$. Thus, $\delta^*(q)=0$ holds
  on $U$ if and only the $(x,y)$--parameter lines are curvature lines of $f$.
  Furthermore, it follows that the zeros of $q$ are umbilics of $f$.
\end{proof}

\begin{Rem}
  The proof of Proposition~\ref{pro:isothermic_is_deltastar_of_q_eq_0}
  shows that a strongly isothermic immersion is isothermic in the
  sense that, away from umbilics, it admits conformal curvature line
  coordinates.
\end{Rem}

Examples of strongly isothermic immersions are surfaces of revolution and
constant mean curvature surfaces in space forms. For the latter, a
holomorphic quadratic differential satisfying $\delta^*(q)=0$ is given by
$q=iQ$, where $Q$ denotes the Hopf differential.

\begin{Pro}\label{lem:strongly_isothermic}
  Strongly isothermic immersions of a compact surface are those points in the
  space of immersions at which the projection to Teichm\"uller space does not
  have maximal rank.
\end{Pro}

\begin{proof}
  Let $M$ be a compact, oriented surface of genus $g\geq1$.  The Teich\-m\"uller
  space $\T(M)$ of $M$ is the quotient of the set of complex structures
  \[ \mathcal C(M)=\{\,J\in\Gamma(\End(TM)) \mid J^2=-\Id\,\}\] on $M$ by the
  group of diffeomorphisms of $M$ that are homotopic to the identity.
  It is well known, see e.g.~\cite{EaEe69,Tr92}, that $\T(M)$ is a
  finite dimensional smooth manifold and that the projection
  $\Phi\colon\mathcal C(M)\to \T(M)$ is differentiable. The kernel of its
  differential at $J$ is the image of $X\in \Gamma(TM)\mapsto \mathcal L_XJ=\dbar
  X$, i.e., the space that, with respect to the pairing
  \eqref{eq:pairing2}, is perpendicular to the finite dimensional
  vector space of holomorphic quadratic differentials of the Riemann
  surface $(M,J)$:
  \begin{equation}\tag{A.1}\label{eq:ker_dPhi}
    \ker(d\Phi_J)=\dbar(\Gamma(TM))=H^0(K^2)^\perp.
  \end{equation}

  An immersion $f\colon M\to \bar M$ induces a unique complex structure
  $J\in\mathcal C(M)$ compatible with the orientation and the induced metric.
  Let $\tau$ denote the composition of the map $f\mapsto J$ with $\Phi$.  By
  Lemma~\ref{lem:change_of_J} and \eqref{eq:delta}, for every infinitesimal
  variation $\dot f = u \xi + df(X)$ of $f$, the differential $d\tau_f$ of
  $\tau$ at $f$ is
  \[ d\tau_f (\dot f)=d\Phi_J\circ \delta(u).\] Hence, $d\tau_f$ is surjective if and
  only if $\ker(d\Phi_J)+\Im(\delta) = \Gamma(\End_-(TM))$, which, by
  \eqref{eq:ker_dPhi}, is  equivalent to $\ker \delta^*\cap H^0(K^2)=0$, because
  $H^0(K^2)$ has finite dimension.
\end{proof}

\begin{The}\label{the:appendix}
  Let $f\colon M\to \bar M$ be an immersion of a compact Riemann surface. If
  $f$ is not strongly isothermic, then every infinitesimal conformal variation of
  $f$ is the derivative of a genuine conformal variation.
\end{The}

As a direct consequence we obtain the following corollary.

\begin{Cor}\label{cor:non_isothermic}
  An immersion $f\colon M\to \bar M$ of a compact Riemann surface that is not
  strongly isothermic is constrained $\Func$--critical if and only if it is a
  critical point of $\Func$ under all genuine conformal variations.
\end{Cor}

\begin{proof}[Proof of Theorem~\ref{the:appendix}]
  By Proposition~\ref{lem:strongly_isothermic}, if $f$ is not strongly isothermic,
  the differential $d\tau_{f}$ of the map $\tau$ at $f$ is surjective. Hence,
  there exist functions $u_1,\ldots, u_n\in C^\infty(M)$ such that the images
  of the corresponding infinitesimal normal variations $d\tau_{f}
  (u_1\xi),\ldots, d\tau_{f}(u_n\xi )$ span the tangent space to
  Teich\-m\"uller space at $\tau(f)$. Let $\dot f = u_0 \xi + df(X)$ be an
  infinitesimal conformal variation of $f$, i.e., $d\tau_{f}(\dot
  f)=d\Phi\circ \delta(u_0)=0$.  Then, in a neighborhood of zero in $\R^{n+1}$
  the maps
  \begin{equation*}
    f_{t_0,\ldots,t_n}(p)=\exp_{f(p)}(t_0(u_0\xi +df(X))+(t_1u_1\ldots+t_nu_n)\xi)
  \end{equation*}
  are immersions and $\tilde \tau (t_0,\ldots,t_n)=\tau(f_{t_0,\ldots,t_n})$
  is a submersion onto Teich\-m\"uller space. The implicit function theorem
  implies that there exists a curve \mbox{$\gamma\colon]-\epsilon,\epsilon[\to \R^n$} such that
  $\gamma(0)=0$ and $\tilde\tau(t,\gamma(t)) =\tilde\tau (0)=\tau(f)$. In
  particular, taking the derivative of $\tilde\tau (t,\gamma(t))$ at $t=0$
  yields $\sum_{i=1}^n\dot \gamma_id\tau_{f} (u_i\xi)=0$ which implies $\dot
  \gamma_i=0$.  Hence, the genuine conformal variation
  $f_t=f_{t,\gamma_1(t)\ldots,\gamma_n(t)}$ of $f$ satisfies $\dot f={\frac d{dt}}_{|t=0}
  f_t=u_0\xi+df(X)$.
\end{proof}

\subsection*{Acknowledgments}

We would like to thank Fran Burstall and Franz Pedit for helpful discussions
and Friederike Sziegoleit for the pictures of Figure 2.


\end{document}